\documentclass[graybox]{svmult}

\usepackage{type1cm}        
\usepackage{makeidx}         
\usepackage{graphicx}        
\usepackage{multicol}        
\usepackage[bottom]{footmisc}

\usepackage{newtxtext}       %
\usepackage{newtxmath}       

\usepackage{algorithm}
\usepackage{algpseudocode}
\usepackage{subfigure}
\algnewcommand{\LIf}[1]{\State\algorithmicif\ {\footnotesize #1}\ \algorithmicthen}
\algnewcommand{\EndLIf}{\unskip\ \algorithmicend\ \algorithmicif}

\def \mypara {\ensuremath{m}}
\newcommand{\mystate}[1][\myhyper]{\ensuremath{u_{#1}}}
\def \mydata {\ensuremath{d}}
\def \mynoise {\ensuremath{\eta }}
\def \myRBstate {\ensuremath{u_{\myhyper, \rm{R}}}}
\def \mygenstateRB {\ensuremath{u_{\rm{R}}}}
\def \mygenstate {\ensuremath{u}}

\def \myhyper {\ensuremath{ \theta}}
\def \myhypereps {\ensuremath{ \varepsilon_{\myhyper} }}
\def \mylambda {\ensuremath { \sigma^2}}
\def \mylambdasqrt {\ensuremath {\sigma}}

\def \mytestfct {\ensuremath{\psi}}
\newcommand{\myeigval}[1][i]{\ensuremath{\lambda_{#1}}}
\newcommand{\myeigvec}[1][i]{\ensuremath{ \mypara_{\myeigval[#1]} }}

\newcommand{\myetainf}[1][\myparaspace]{\ensuremath{\underline{\eta}{}_{#1,\myhyper}}}
\newcommand{\myetasup}[1][\myparaspace]{\ensuremath{\overline{\eta}_{#1,\myhyper}}}

\def \mybetaT {\ensuremath {\beta_{\myhyper, \myobs}}}
\newcommand{\mybetaTtwo}[2]{\tilde{\beta}_{\myobs, #1, #2}}
\def \mybetaTRB {\ensuremath {\beta_{\myhyper, \myobs, \rm{R}}}}

\def \myCone {\ensuremath {C_{\myhyper, \myobs, \mylambda} }}

\def \mygammaobs {\ensuremath { \gamma_{\myobs} }}
\def \myratio {\ensuremath {C_{\mypriorCovinv}}}

\def \mydatadim {\ensuremath{K}}
\def \mylibrarydim {\ensuremath{\mydatadim_{\mylibrary}}}
\def \myhyperdim {\ensuremath{p}}
\def \myparadim {\ensuremath {M}}

\def \myparaspace{\ensuremath{\mathbb{R}^{\myparadim}}}
\def \myparasubspace{\ensuremath{X}}
\def \myparaspacehyper{\ensuremath{ X_{\myhyper}}}
\def \myparaspacehyperortho{\ensuremath{ X_{\myhyper}^{\perp}}}

\def \mystatespace {\ensuremath{ \mathcal{U}}}
\def \mystatespacedual {\ensuremath{ \mystatespace ' }}

\def \mydataspace{\ensuremath{ \mathbb{R}^{\mydatadim} }}
\def \myhyperdomain{\ensuremath{ \mathcal{P}}}
\def \myhyperdomainspace{\ensuremath{ \mathbb{R}^{\myhyperdim}}}

\def \mylibrary{\ensuremath{ \mathcal{L} }}
\def \myXitrain{\ensuremath{ \Xi_{\rm{train}}}}
\def \mytargetbeta{\ensuremath{ \beta_0 }}
\def \myiterbeta{\ensuremath{ \beta }}
\newcommand{\mybetaTRBhyper}[1]{\ensuremath {\beta_{\myobs, \rm{R}}(#1)}}
\def \mydatadimmax{\ensuremath{ \mydatadim_{\rm{max}} }}
\def \mydatadimiter{\ensuremath{ K }}
\newcommand{\myRBstatehyper}[1]{\ensuremath{u_{#1, \rm{R}}}}
\def \mystateiter{\ensuremath{ u_{\mydatadimiter+1} }}

\newcommand{\mystateinner}[2]{{(#1, #2)_{\mystatespace}}}
\newcommand{\mystatenorm}[1]{|| #1 ||_{\mystatespace}}

\newcommand{\mypriorinner}[2]{(#1, #2)_{\mypriorCovinv}}
\newcommand{\mypriornorm}[1]{|| #1 ||_{\mypriorCovinv}}
\newcommand{\mynoisenorm}[1]{|| #1 ||_{\mynoiseCovinvobs}}
\newcommand{\mynoiseinner}[2]{{(#1, #2)_{\mynoiseCovinvobs}}}

\newcommand{\myMnorm}[1]{|| #1 ||_{\mathbb{R}^{\myparadim}}}

\def \mylhs {{\ensuremath a_{\myhyper}}}
\def \myrhs {{\ensuremath b_{\myhyper}}}

\def \myG {\ensuremath {G_ {\myhyper, \myobs}}}
\def \myGadj {\ensuremath {\myG^*}}
\def \myobs {\ensuremath {L}}
\newcommand{\myobsl}[1]{\ensuremath {l_{#1}}}
\newcommand{\mylibraryl}[1][i]{\ensuremath {l_{#1}}}
\def \mygenMat {\ensuremath {}}

\def \mygenProj {\ensuremath {\Pi_{\myparaspacehyperortho}}}

\def \mypriormeasure {\ensuremath {\mu_{0}}}
\def \mypriorpara {\ensuremath{ \mypara_{0}}}
\def \mypriorCov {\ensuremath {\Sigma_0}}
\def \mypriorCovinv {\ensuremath {\mypriorCov^{-1}}}

\def \mynoiseCov {\ensuremath {\Sigma_{\rm{noise}}}}
\def \mynoiseCovobs {\ensuremath {\Sigma_{\myobs}}}

\def \mynoiseCovinvobs {\ensuremath {\Sigma_{\myobs}^{-1}}}

\def \myposterior {\ensuremath {\pi_{\rm{post}}}}
\def \myposteriormeasure {\ensuremath {\mu^{\myobs,\mydata}}}
\def \myposteriorpara {\ensuremath {\mypara_{\rm{post}}^{\myhyper, \myobs}}}
\def \myposteriorCov {\ensuremath {\Sigma_{\rm{post}}^{\theta, \myobs}}}

\def \myOmega {\ensuremath {\Omega}}
\def \myGammain {\ensuremath {\Gamma_{\rm{in}}}}
\def \mypoly {\ensuremath {p_i}}

\usepackage[colorinlistoftodos,prependcaption,textsize=tiny]{todonotes}
\usepackage{soul,xargs}

\newcommandx{\missing}[2][1=]{\todo[linecolor=yellow,backgroundcolor=yellow!25,bordercolor=yellow,#1]{#2}}
\newcommandx{\unsure}[2][1=]{\todo[linecolor=red,backgroundcolor=red!25,bordercolor=red,#1]{#2}}
\newcommandx{\change}[2][1=]{\todo[linecolor=blue,backgroundcolor=blue!25,bordercolor=blue,#1]{#2}}
\newcommandx{\info}[2][1=]{\todo[linecolor=green,backgroundcolor=green!25,bordercolor=green,#1]{#2}}
\newcommandx{\improvement}[2][1=]{\todo[linecolor=violet,backgroundcolor=violet!25,bordercolor=violet,#1]{#2}}
\newcommandx{\thiswillnotshow}[2][1=]{\todo[disable,#1]{#2}}

\makeindex             

\begin{document}

\title*{A sequential sensor selection strategy for hyper-parameterized linear Bayesian inverse problems}
\titlerunning{Sensor selection strategy for hyper-parameterized linear Bayesian inverse problems}
\author{Nicole Aretz-Nellesen, Peng Chen, Martin A. Grepl and Karen Veroy}
\institute{Nicole Aretz-Nellesen \at International Research Training Group Modern Inverse Problems, RWTH Aachen University, Schinkelstraße 2, 52062 Aachen, Germany, \email{nellesen@aices.rwth-aachen.de}
\and Peng Chen \at Oden Institute of Computational Engineering Science, UT Austin, 201 E 24th St, Austin, TX 78712, USA, \email{peng@ices.utexas.edu}
\and Martin A. Grepl \at Numerical Mathematics (IGPM), RWTH Aachen University, Templergraben 55, 52056 Aachen,
Germany, \email{grepl@igpm.rwth-aachen.de}
\and Karen Veroy \at Aachen Institute for Advanced Study in Computational Engineering Science (AICES), RWTH Aachen University, Schinkelstraße 2, 52062 Aachen, Germany, \email{veroy@aices.rwth-aachen.de}}
%
%
\maketitle

\abstract{
We consider optimal sensor placement for hyper-parameterized linear Bayesian inverse problems, where the hyper-parameter characterizes nonlinear flexibilities in the forward model, and is considered for a range of possible values.
This model variability needs to be taken into account for the experimental design to guarantee that the Bayesian inverse solution is uniformly informative.
In this work we link the numerical stability of the maximum a posterior point and A-optimal experimental design 
to an observability coefficient that directly describes the influence of the chosen sensors.
We propose an algorithm that iteratively chooses the sensor locations to improve this coefficient and thereby decrease the eigenvalues of the posterior covariance matrix.
This algorithm exploits the structure of the solution manifold in the hyper-parameter domain via a reduced basis surrogate solution for computational efficiency.
We illustrate our results with a steady-state thermal conduction problem.
}

\section{Introduction}
\label{sec:Introduction}
Mathematical models of physical processes often depend on parameters, such as material properties or source terms, that are known only with some uncertainty.
Experimental measurement data can help estimate these parameters and thereby improve the meaningfulness of the model.
The Bayesian approach to inverse problems (cf. \cite{Stuart_2010a}) yields a (posterior) probability distribution for these parameters that reflects both the prior distribution in the parameters and measurement data.

A major challenge in inverse problems is sensor placement to obtain informative measurement data at restricted experimental cost.
There exists a vast optimal experimental design (OED) community focused on different problem types and optimal design criteria.
The literature most related to this contribution is the discussion of A-optimality for infinite-dimensional linear Bayesian inverse problems in \cite{Alexanderian_2014a, Alexanderian_2016a}, and the greedy orthogonal matching pursuit algorithm for data assimilation in \cite{Binev_2018a, Maday_2015a}.

In this paper, we consider the optimal placement of sensors to infer a parameter from noisy data in a linear Bayesian inverse problem subject to flexible hyper-parameters.
The hyper-parameters characterize variability of the forward model, e.g. variable material properties or geometry, that needs to be taken into account for the sensor placement.
For instance, the Bayesian inference problem might need to be solved for multiple data sets with known hyper-parameters, or a most suitable hyper-parameter might need to be sought for fixed data in an "outer loop" optimization.
In either case, the same sensors are used within all inference problems, thus necessitating a uniformly "good" choice.
The objective of this paper is to provide a sensor selection strategy that follows, uniformly for all hyper-parameters, the A-optimal design criterion of minimizing the trace of the posterior covariance matrix.

In \cite{Aretz_2019a}, we developed and utilized a numerical stability analysis for parameterized 3D-VAR data assimilation over a linear model correction term to find design criteria for stability-based sensor selection.
In this contribution, we first re-interpret these results in the hyper-parameterized linear Bayesian inversion setting, and then show their relation to A-optimal experimental design.
This analysis leads to a greedy algorithm that iteratively chooses sensor locations that, under certain assumptions, uniformly decrease the trace of the posterior covariance matrix.

In the upcoming section, we specify our linear forward model, and pose the
Bayesian inversion problem in a hyper-parameterized context. 
We then, in Sec. \ref{sec:StabilityOED}, show the link between its numerical stability, different model coefficients, the eigenvalues of the posterior covariance matrix, and A-optimal experimental design.
In Sec. \ref{sec:sensors} we propose an algorithm to exploit this connection, and present numerical results in Sec. \ref{sec:implementation} for a thermal conduction problem.
We conclude in Sec. \ref{sec:conclusion}.

\section{A Hyper-Parameterized Bayesian Inverse Problem}
\label{sec:setting}

We consider a linear Bayesian inverse problem setting for the inference of a finite-dimensional\footnote{
The extension to the infinite-dimensional setting poses additional challenges that will be discussed in a future work.
}
parameter\footnote{
A general finite-dimensional space can be considered via an affine transformation, c.f. \cite{Alexanderian_2016a, DaPrato_2006a}
}
$\mypara \in \myparaspace$ from noisy data $\mydata \in \mydataspace$ subject to different hyper-parameters $\myhyper$ that characterize nonlinear (in $\myhyper$) flexibility in the linear (in $\mypara$) forward model.
Our objective is to find conditions for an observation operator that is uniformly informative for all hyper-parameters.
In the following, we specify the forward model and the Bayesian inverse problem, before analysing it in Sec. \ref{sec:StabilityOED}.

Following the Bayesian approach to inverse problems, we consider $\mypara$ to be a random variable, and model our prior belief in its distribution through a non-degenerate Gaussian prior measure $\mypriormeasure = \mathcal{N}(\mypriorpara, \mypriorCov)$ with mean $\mypriorpara$ and symmetric positive-definite (s.p.d.) covariance $\mypriorCov \in \mathbb{R}^{\myparadim \times \myparadim}$.
We define the inner product $\mypriorinner{\mypara_1}{\mypara_2} := \mypara_1^T \mypriorCovinv \mypara_2$ and norm $\mypriornorm{\mypara_1}^2:= \mypriorinner{\mypara_1}{\mypara_1}$
for $\mypara_1, \mypara_2 \in \myparaspace$.

For the forward model, let $(\mystatespace, \mystateinner{\cdot}{\cdot})$ be a Hilbert space with induced norm $\mystatenorm{\mygenstate}^2 := \mystateinner{\mygenstate}{\mygenstate}$, and let $\myhyperdomain \subset \myhyperdomainspace$ be a compact set of possible hyper-parameters.
For any $\myhyper \in \myhyperdomain$, we let
$\mylhs : \mystatespace \times \mystatespace \rightarrow \mathbb{R}$ and
$\myrhs : \myparaspace \times \mystatespace \rightarrow \mathbb{R}$ be non-trivial bilinear forms that are affine\footnote{For conciseness, we refer the reader to \cite{Aretz_2019a} for a definition of these properties.}
 and bounded uniformly in $\myhyper$, with the additional assumption that $\mylhs$ is also uniformly coercive.\footnote{
We can readily generalize this setting to non-coercive problems by employing a Petrov-Galerkin formulation. A stability analysis similar to \cite{Aretz_2019a} will be explored in a future publication.}
Under these assumptions there exists, for any parameter $\mypara \in \myparaspace$, a unique, bounded solution to the problem
\begin{equation}\label{eq:model}
\text{find } \mystate(\mypara) \in \mystatespace \quad \text{ s.t. } \quad
\mylhs(\mystate, \mytestfct) = \myrhs(\mypara, \mytestfct) \quad \forall \mytestfct \in \mystatespace.
\end{equation}
We define, for $\myparasubspace \subset \myparaspace$, the ratios $\myetainf[\myparasubspace] := \inf_{\mypara \in \myparasubspace} \mystatenorm{\mystate(\mypara)}/\mypriornorm{\mypara} \ge 0$ and $\myetasup[\myparasubspace] := \sup _{\mypara \in \myparasubspace} \mystatenorm{\mystate(\mypara)}/\mypriornorm{\mypara} < \infty$.
Moreover, we define the closed subspace
\begin{equation}\label{eq:myparaspacehyper}
\myparaspacehyper := \{ \mypara \in \myparaspace: ~ \mystate(\mypara) = 0\} = \{\mypara \in \myparaspace: ~\myrhs(\mypara, \cdot ) = 0\} \subset \myparaspace
\end{equation}
of all parameter directions that do not change the state, and let $\myparaspacehyperortho$ denote its orthogonal complement in the Euclidean inner product.
In particular, we have $\myetainf[\myparaspacehyperortho] > 0$.

For our sensors, we consider a library $\mylibrary = \{\mylibraryl[k] \}_{k=1}^{\mylibrarydim}$ of $\mylibrarydim < \infty$ sensors $\mylibraryl[k] \in \mystatespacedual$.
For a selection $\mylibraryl[k_1], \dots, \mylibraryl[k_\mydatadim] \in \mylibrary$ of these sensors, we define the observation operator $\myobs = (\mylibraryl[k_1], \cdots, \mylibraryl[k_\mydatadim])^T : \mystatespace \rightarrow \mydataspace$.
Measurement data for a parameter $\mypara \in \myparaspace$ is obtained by applying $\myobs$ to the state $\mystate(\mypara)$.
This gives us the linear, bounded parameter-to-observable map $\myG : \myparaspace \rightarrow \mydataspace$, $\myG(\mypara) := \myobs \mystate(\mypara)$.
Our objective is to choose $\myobs$ from $\mylibrary$ so that it is approximately A-optimal over $\myhyperdomain$.

For the noise model, we assume to be given an s.p.d. covariance matrix $\mynoiseCov \in \mathbb{R}^{\mylibrarydim \times \mylibrarydim}$ that describes how the observation noise between all sensors in $\mylibrary$ is correlated.
The covariance for the sensors in $\myobs$ is then described by the submatrix $\mynoiseCovobs \in \mathbb{R}^{\mydatadim \times \mydatadim}$ with $(\mynoiseCovobs)_{i,j} = (\mynoiseCov)_{k_i, k_j}$.
For fixed $\myobs$ and for data $\mydata_1, \mydata_2 \in \mathbb{R}^{\mydatadim}$ we define the inner product $\mynoiseinner{\mydata_1}{\mydata_2} := \mydata_1^T \mynoiseCovinvobs \mydata_2$ and induced norm $\mynoisenorm{\mydata_1}^2:= \mynoiseinner{\mydata_1}{\mydata_1}$.
We define $\mygammaobs := \sup_{\mygenstate \in \mystatespace} \mynoisenorm{\myobs \mygenstate}/\mystatenorm{\mygenstate}$ as the norm of $\myobs$.
We model the data to be of the form
\begin{equation}\label{eq:noisemodel}
\mydata = \myG (\mypara) + \mynoise
\quad \text{with Gaussian additive noise} \quad
\mynoise \sim \mathcal{N}(0, \mylambda \mynoiseCovobs)
\end{equation}
and scaling parameter $\mylambdasqrt > 0$.
For given data $\mydata \in \mydataspace$ from an observation operator $\myobs$, the posterior probability density function of the posterior measure $\myposteriormeasure$ is then given through Bayes' theorem by
\begin{equation}\label{eq:posterior}
\myposterior (\mypara | \mydata) \propto \exp \big{(} - \textstyle\frac{1}{2\mylambda} \mynoisenorm{\myG(\mypara) - \mydata}^2 - \textstyle \frac{1}{2} \mypriornorm{\mypara - \mypriorpara}^2 \big{)},
\end{equation}
where we omit the normalization constant $\textstyle Z = \int_{\myparaspace} \exp(-\frac{1}{2\mylambda} \mynoisenorm{\myG(\mypara) - \mydata}^2) d\mypriormeasure$.
Since $\textstyle \myG$ is linear, the posterior is a Gaussian (see, e.g., \cite[Thm 2.4]{Stuart_2010a}), $\textstyle \myposteriormeasure = \mathcal{N}(\myposteriorpara(\mydata), \myposteriorCov)$, with mean $\textstyle \myposteriorpara(\mydata) = \myposteriorCov \left(\textstyle \frac{1}{\mylambda} \myGadj \mynoiseCovinvobs \mydata + \mypriorCovinv \mypriorpara \right)$ and covariance matrix $\textstyle \myposteriorCov = \left( \textstyle \frac{1}{\mylambda}\myGadj \mynoiseCovinvobs \myG + \mypriorCovinv \right)^{-1}$.

\section{Numerical Stability and A-Optimal Experimental Design}
\label{sec:StabilityOED}

In the following, we first comment on the connection between the numerical stability of the MAP point and the observation operator $\myobs$.
We then link this analysis to A-optimal experimental design.

Since $\myposteriormeasure$ is Gaussian, its mean $\myposteriorpara$ is the maximum a posteriori (MAP) point, and hence the solution to the minimization problem
\begin{equation}\label{eq:min}
\min _{\mypara \in \myparaspace} \textstyle \frac{1}{2\mylambda} \mynoisenorm{\myobs \mystate (\mypara) - \mydata}^2 + \textstyle \frac{1}{2} \mypriornorm{\mypara - \mypriorpara}^2.
\end{equation}
Through a reformulation as a saddle-point problem, the numerical stability of \eqref{eq:min} can be analyzed with respect to $\myhyper$, $\myobs$, and $\mylambda$ (c.f. \cite{Aretz_2019a} for an analogous analysis).
In particular, the difference in the MAP points and states for different $\mydata_1, \mydata_2 \in \mydataspace$ is bounded by the difference in data.
We have, with $\tilde{\mypara}(\mydata) := \myposteriorpara(\mydata)$ for readability,
\begin{equation}
\mypriornorm{\tilde{\mypara}(\mydata_1)-\tilde{\mypara}(\mydata_2)}^2 + \mystatenorm{\mystate(\tilde{\mypara}(\mydata_1))-\mystate(\tilde{\mypara}(\mydata_2))}^2 \le (\myCone) ^2 \mynoisenorm{\mydata_1 - \mydata_2}^2.
\end{equation}
The stability coefficient $\myCone > 0$ quantifies the influence of noise on the MAP point.
It has the form $\myCone = \mygammaobs (1+\eta^2)/(\mylambda + \mybetaT^2 \eta^2)$, where $\eta = \myetasup$ if $\mybetaT^2 \le \mylambda$, and $\eta = \myetainf$ otherwise, and $\mybetaT$ is the observability coefficient
\begin{equation}\label{eq:betaT}
\mybetaT := \inf \{\mynoisenorm{\myobs \mystate(\mypara)}:~ \mystatenorm{\mystate(\mypara)} = 1, ~ \mypara \in \myparaspace \} .
\end{equation}
$\myCone$ decreases in $\mybetaT$, and remains bounded for $\mylambda \rightarrow 0$ iff $\myetainf > 0$ and $\mybetaT > 0$.
Increasing $\mybetaT$ can hence help improve robustness of $\myposteriorpara$ against noise.

The goal in A-optimal experimental design is to choose sensors to minimize the trace of the posterior covariance matrix $\myposteriorCov$.
Geometrically, this corresponds to minimizing the mean axis of the uncertainty ellipsoid (c.f. \cite{Ucinski_2004a}).
In the following, we bound the eigenvalues of $\myposteriorCov$ via $\mybetaT$ and $\myetainf[\myparaspacehyperortho]$.
These can then be used to choose sensors to decrease the bounds of the eigenvalues, and consequently of the trace.

Utilizing the definitions of $\myG$ and $\myetainf[\myeigvec]$ in Sec. \ref{sec:setting}, and $\mybetaT$ in \eqref{eq:betaT}, we observe
\begin{equation} \label{eq:updatebound}
\begin{aligned}
\mypara^T\myGadj \mynoiseCovinvobs \myG \mypara 
&= \mynoisenorm{ \myobs \mystate(\mypara)}^2
\ge \mybetaT^2 \mystatenorm{\mystate(\mypara)}^2
= \mybetaT^2 \mystatenorm{\mystate(\mygenProj \mypara)}^2\\
&\ge \mybetaT^2 \myetainf[\myparaspacehyperortho]^2 \mypriornorm{\mygenProj \mypara}^2 
\ge \mybetaT^2 \myetainf[\myparaspacehyperortho]^2 \myratio^{-2} \myMnorm{\mygenProj \mypara}^2, 
\end{aligned}
\end{equation}
where $\mygenProj$ is the orthogonal projection onto $\myparaspacehyperortho$ in the Euclidean inner product and $\textstyle \myratio := \sup _{\mypara \in \myparaspace} \myMnorm{\mypara}/\mypriornorm{\mypara}$ is a norm equivalence constant.

Let $0 < \myeigval[1] \le \dots \le \myeigval[\myparadim]$ be the eigenvalues of the posterior covariance matrix $\textstyle \myposteriorCov$, including duplicates.
Since $\textstyle \myposteriorCov$ is s.p.d., there exists an orthonormal eigenvector basis $\textstyle (\myeigvec)_{i=1}^{\myparadim}$ of $\textstyle \myparaspace$, i.e. $\textstyle{ m_{\lambda_i}^T \myeigvec[j] = \delta_{i,j} }$ and $\textstyle{ \myposteriorCov \myeigvec = \myeigval \mygenMat \myeigvec}$.
With the explicit formula for $\myposteriorCov$ from Sec. \ref{sec:setting}, the last equation is equivalent to $\textstyle \frac{1}{\myeigval} \mygenMat \myeigvec = \textstyle \frac{1}{\mylambda}\myGadj \mynoiseCovinvobs \myG \myeigvec+ \mypriorCovinv \myeigvec $.
Premultiplying by $\myeigvec^T$ and inserting \eqref{eq:updatebound} yields
$\frac{1}{\myeigval} = \frac{1}{\myeigval} \myMnorm{\myeigvec}^2 \ge (\textstyle \frac{1}{\mylambda}\mybetaT^2 \myetainf[\myparaspacehyperortho]^2 \myratio^{-2} \myMnorm{\mygenProj \myeigvec}^2  + \myratio^{-2})$, and hence 
\begin{align*}
\textstyle \myeigval \le \myratio^2/(\textstyle \frac{1}{\mylambda} \mybetaT^2 \myetainf[\myparaspacehyperortho]^2 \myMnorm{\mygenProj \myeigvec}^2 + 1).
\end{align*}
Summing over all eigenvalues, including duplicates, we can now bound 
\begin{align*}
\text{trace}(\myposteriorCov) = \textstyle{\sum_{i=1}^{\myparadim}} \myeigval \le \myratio^2 \sum_{i=1}^{\myparadim} (\textstyle \frac{1}{\mylambda} \mybetaT^2 \myetainf[\myparaspacehyperortho]^2 \myMnorm{\mygenProj \myeigvec}^2 + 1)^{-1}.
\end{align*}
Although $\myMnorm{\mygenProj \myeigvec} \ge 0$ is unknown for any individual $\myeigval$, by exploiting that $(\myeigvec)_{i=1}^{\myparadim}$ is an orthonormal basis, it can be shown that $\sum _{i=1}^{\myparadim} \myMnorm{\mygenProj \myeigvec}^2 = \dim \myparaspacehyperortho$.
Our strategy is to choose $\myobs$ to increase $\mybetaT$; this decreases the bound for each $\myeigval$ with $\myMnorm{\mygenProj \myeigvec} > 0$, and hence also the bound for the trace.
The coefficient $\mybetaT$ becomes more influential the more the data is trusted, i.e. for $\mylambda$ small.

\section{Sensor Selection Strategy}
\label{sec:sensors}

Our goal is to choose sensors $\{\myobsl{k}\}_{k=1}^{\mydatadim}$ so that $\mybetaT$ is uniformly large over the hyper-parameter domain $\myhyperdomain$.
The major challenge for achieving this goal is that evaluating $\mybetaT$ for any $\myhyper$ involves solving the forward problem \eqref{eq:model} for each basis vector of $\myparaspace$.
We address this problem by approximating the solution $\mystate(\mypara)$ of \eqref{eq:model} with a surrogate reduced basis (RB) solution $\myRBstate(\mypara)$, that can be computed at a considerably reduced computational cost for a specified accuracy:
Suppose for every hyper-parameter $\myhyper \in \myhyperdomain$ and every parameter $\mypara \in \myparaspace$, we can compute $\myRBstate(\mypara) \in \mystatespace$ such that $\mystatenorm{\mystate(\mypara)-\myRBstate(\mypara)} \le \myhypereps \mystatenorm{\mystate(\mypara)}$ for a relative accuracy $0 \le \myhypereps \le \varepsilon < 1$.
It can then be shown analogously to \cite[sec. 5.1]{Aretz_2019a} that $\mybetaT \ge (1-\myhypereps) \mybetaTRB - \mygammaobs \myhypereps$, where $\mybetaTRB$ is defined over the surrogate model analogously to \eqref{eq:betaT}, and $\mygammaobs$ is the norm of $\myobs$.
The upper bound $\varepsilon < 1$ ensures that $\myRBstate(\mypara) = 0$ iff $\mystate(\mypara) = 0$; therefore $\mybetaTRB$ is defined over the same parameter subspace $\myparaspacehyperortho$ as $\mybetaT$.
Supposing $\varepsilon$ is small enough, we propose to exploit the lower bound of $\mybetaT$ by choosing the sensors in $\myobs$ via an iterative greedy approach over $\myhyperdomain$ to increase $\mybetaTRB$, and subsequently $\mybetaT$.
\begin{algorithm}[h]
\caption{Stability-based sensor selection}
\label{alg:greedyOMP}
\textbf{Given:} a training set $\myXitrain$, a library $\mylibrary$, a target value $\mytargetbeta$, a starting parameter $\myhyper_1 \in \myXitrain$, and an upper limit $\mydatadimmax$ to the number of sensors.
\begin{algorithmic}[1]
\State $\myobs \gets \{0\}$, $\myiterbeta \gets 0$, $\mydatadimiter \gets 0$
\While {$\myiterbeta < \mytargetbeta$ and $\mydatadimiter < \mydatadimmax$}
\State $\mystateiter \gets \text{arg} \min \{ \mynoisenorm{\myobs \mygenstateRB}/\mystatenorm{\mygenstateRB}: ~ \mygenstateRB = \myRBstatehyper{\myhyper_{\mydatadimiter+1}}(\mypara) \text{ for } \mypara \in \myparaspace\}$  \label{line:choosestate}
\State choose $l \in \mylibrary$ such that 
$||[\myobs, l] \mystateiter||_{\Sigma_{[\myobs, l]}^{-1}}$ is maximal \label{line:choosesensor}
\State $\myobs \gets [\myobs, l]$, $\mydatadimiter \gets \mydatadimiter + 1$
\State $\myhyper_{\mydatadimiter+1} \gets \text{arg} \min _{\myhyper \in \myXitrain} \mybetaTRB$, $\myiterbeta \gets \mybetaTRBhyper{\myhyper_{\mydatadimiter + 1}}$ \label{line:next}
\EndWhile
\end{algorithmic}
\end{algorithm}

Following the ideas in \cite{Binev_2018a, Maday_2015a}, in each iteration of the loop, the algorithm first chooses (line \eqref{line:choosestate}) the state $\mystateiter$ which realizes the minimum observability coefficient; it then searches the library $\mylibrary$ for the best sensor to observe this state (line \eqref{line:choosesensor}), and then extend the observation operator.
Line \eqref{line:choosesensor} involves first computing $l(\mystateiter)$ (in FE dimension) and then $||[\myobs, l] \mystateiter||_{\Sigma_{[\myobs, l]}^{-1}}$ (in $\mathcal{O}(\mydatadim^3)$) for each $l \in \mylibrary$.
In line \eqref{line:next} the algorithm then finds the hyper-parameter by iterating over $\myhyper \in \myXitrain$ and computing $\mybetaTRB$ via an eigenvalue problem.
Any computation of $\mybetaTRB$ involves solving the RB problem for each basis vector of $\myparaspace$.
The algorithm terminates when either a maximum number of sensors or a target value \footnote{
Possibilities for target values are highly dependent on the library $\mylibrary$.
In practice, $\mytargetbeta$ should be chosen by carefully monitoring the changes in $\myiterbeta$.
} $\mytargetbeta$ has been reached.

The uniform increase of $\mybetaTRB$ over $\myhyperdomain$ relies on the property that extending the observation operator with a sensor does not decrease $\mybetaTRB$ at different hyper-parameters.
This property is straightforward to prove for uncorrelated noise, but more involved for the general case. 
We will explore this aspect in a future publication.

\section{Numerical Results}
\label{sec:implementation}
\begin{figure}[t]
\vspace{-1ex}
\centering
\hspace{-5ex}
\subfigure[Thermal block]{\label{fig:thermalblock} \includegraphics[width = 0.48\textwidth]{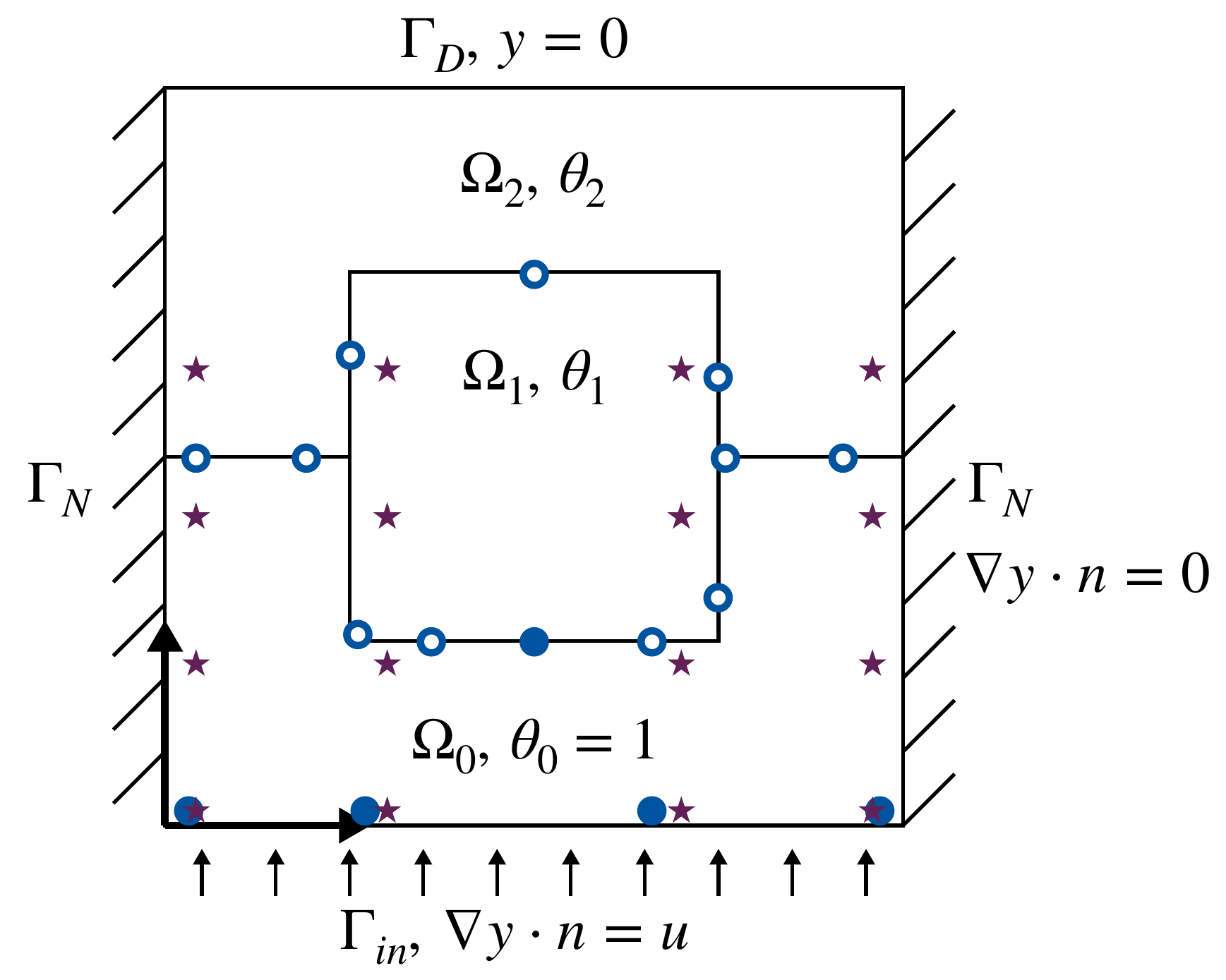} } ~
\subfigure[$\text{trace}(\myposteriorCov)$ vs. $\mybetaT$]{\label{fig:betaT_vs_trace_mean} \hspace{-3ex}\includegraphics[width = 0.53\textwidth]{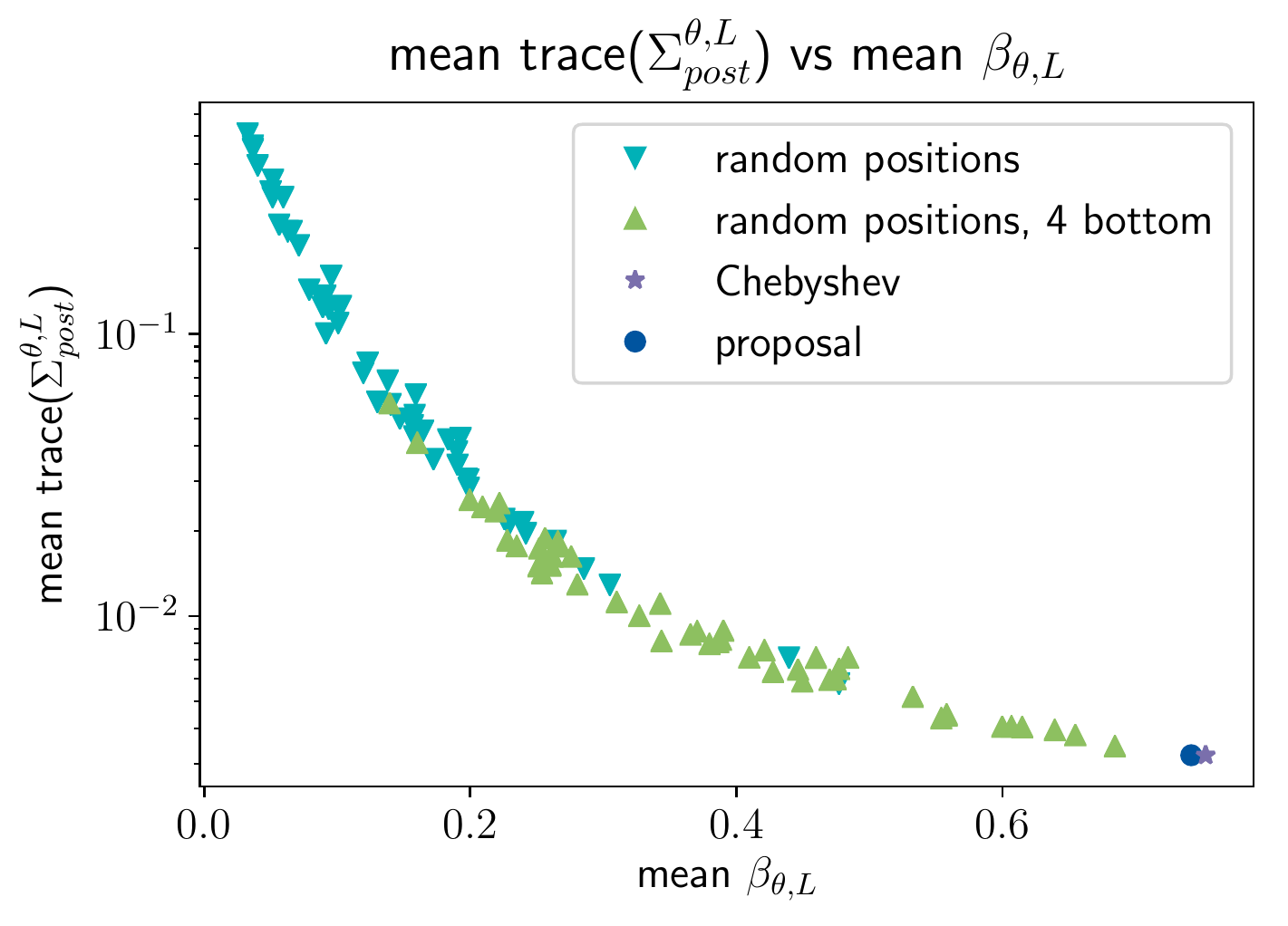} }\hspace{-5ex}
\caption{Fig. \ref{fig:thermalblock}: Domain decomposition of the thermal block problem with boundary conditions.
Sensor centres chosen by Alg. \ref{alg:greedyOMP} are indicated as circles (filled out for $\mybetaT$-criterion), reference Chebyshev centres are marked with stars.
Fig. \ref{fig:betaT_vs_trace_mean}: 
Mean of the trace of the posterior covariance matrix vs. the mean of $\mybetaT$. Mean values are computed over the 1681 hyperparameters in $\Xi_{\rm{test}}$.
}
\end{figure}
We consider a steady-state heat conduction problem $-\myhyper \Delta \mygenstate = 0$ over the unit square $\myOmega$.
The hyper-parameters $\myhyper_i \in [0.1, 10]$ specify different thermal conductivities on three subdomains (overview in Fig. \ref{fig:thermalblock}).
We impose an uncertain inflow boundary condition $\nabla \mygenstate = \sum_{i=0}^3 \mypara_i \mypoly$ on $\myGammain$, where $\mypoly$ is the Legendre polynomial of degree $i$, and $\mypara$ is distributed as $\mathcal{N}((1,0,0,0)^T, I)$ in $\myparaspace = \mathbb{R}^4$, with identity matrix $I$.
We refer to \cite{Aretz_2019a} for a full description of the model problem and algorithm implementation.

Our library consists of functionals $\myobsl{k}(u):= \int _{\myOmega} g_k(x)u(x) dx$, where the $g_k$ are Gaussian functions with standard deviation $0.01$ and centres in a $97 \times 97$ regular grid on $[0.02, 0.98]^2$.
We model $\mynoiseCov$ as the $\mystatespace$-inner-product of the sensors' Riesz representations, and choose $\sigma = 0.01$.
The noise correlation at different sensors is then higher the closer they are placed to each other.
We apply Alg. \ref{alg:greedyOMP} with target value $\mytargetbeta := 0.5$.
In addition to $\mybetaT$, we also consider $
\mybetaTtwo{\myhyper_1}{\myhyper_2} := \inf \{ \mynoisenorm{\myobs \mygenstate}/\mystatenorm{\mygenstate}:~\mygenstate = \mystate[\myhyper_1](\mypara_1) + \mystate[\myhyper_2](\mypara_2), ~ \mypara_1, \mypara_2 \in \myparaspace \}$,
which distinguishes between different hyper-parameters, and selects sensors for "outer-loop" hyper-parameter estimation (c.f. discussion in \cite{Aretz_2019a}).
The 16 selected sensors are indicated in Fig. \ref{fig:thermalblock} as circles.

We compute the mean trace of the posterior covariance matrix and the mean of $\mybetaT$ over a testing set $\Xi_{\rm{test}}$ of $41\times 41$ hyper-parameters located in a regular grid on the logarithmic plane of $\myhyperdomain = (0.1, 10)^2$.
For comparison, we repeat the process 50 times with 16 randomly chosen positions, and another 50 times with 16 randomly chosen positions of which at least 4 are placed close to the inflow boundary with $x_2 = 0.02$.
Also, we consider another set of 16 centers (indicated as stars in Fig. \ref{fig:thermalblock}), where the $x_1$-locations are chosen in our library closest to the Chebyshev interpolation points for polynomials with degree smaller or equal to 3. 
These would be the theoretically optimal points (in $x_1$-direction) for interpolating the Neumann flux.
Fig. \ref{fig:betaT_vs_trace_mean} shows the mean trace and the mean $\mybetaT$ over $\Xi_{\rm{test}}$ for the different sensor sets.

We observe that the sensors chosen by Alg. \ref{alg:greedyOMP} near the inflow boundary are close to the Chebyshev postitions.
Both sensor sets have a similarly high mean value for $\mybetaT$ and a similarly small trace of the posterior covariance.
In contrast, the randomly chosen sensor sets have a larger mean trace and smaller mean $\mybetaT$.
Here the sets chosen with 4 centres near the inflow boundary outperform the completely random ones.
Altogether, we observe a strong correlation between $\mybetaT$ and $\text{trace}(\myposteriorCov)$.

\section{Conclusion}
\label{sec:conclusion}
In this paper we considered a hyper-parameterized linear Bayesian inverse problem and linked its numerical stability analysis to A-optimal experimental design. This analysis permits the development of an algorithm that iteratively chooses sensor locations from a library to reduce the eigenvalues of the posterior covariance matrix uniformly over the hyper-parameter domain. Future work will extend this analysis to Petrov-Galerkin and time-dependent formulations, as well as application to high-dimensional parameter spaces and nonlinear problems.

\begin{acknowledgement}
We would like to thank Tan Bui-Thanh, Youssef Marzouk, Francesco Silva, Andrew Stuart, Dariusz Ucinski, and Keyi Wu for very helpful discussions. This work was supported by the Excellence Initiative of the German federal and state governments and the German Research Foundation through Grants GSC 111 and 33849990/GRK2379 (IRTG Modern Inverse Problems).
\end{acknowledgement}

\bibliographystyle{plain}
\bibliography{literature}
\end{document}